\documentclass[11pt]{article} 
\newcommand{\be}{\begin{equation}} 
\newcommand{\ee}{\end{equation}} 
\newcommand{\bea}{\begin{eqnarray}} 
\newcommand{\eea}{\end{eqnarray}} 
\newcommand{\bean}{\begin{eqnarray*}} 
\newcommand{\eean}{\end{eqnarray*}} 
\newcommand{\brray}{\begin{array}} 
\newcommand{\erray}{\end{array}} 
\newcommand{\ben}{\begin{equation}{nonumber}} 
\newcommand{\een}{\end{equation}{nonumber}} 
\newcommand{\newsection}[1]{\setcounter{dfn}{0} 
\section{#1}}

\newtheorem{dfn}{Definition}[section] 
\newtheorem{thm}[dfn]{Theorem} 
\newtheorem{lmma}[dfn]{Lemma} 
\newtheorem{ppsn}[dfn]{Proposition} 
\newtheorem{crlre}[dfn]{Corollary} 
\newtheorem{xmpl}[dfn]{Example} 
\newtheorem{rmrk}[dfn]{Remark}

\newcommand{\bdfn}{\begin{dfn}} 
\newcommand{\bthm}{\begin{thm}} 
\newcommand{\blmma}{\begin{lmma}} 
\newcommand{\bppsn}{\begin{ppsn}} 
\newcommand{\bcrlre}{\begin{crlre}} 
\newcommand{\bxmpl}{\begin{xmpl}} 
\newcommand{\brmrk}{\begin{rmrk}} 
 
\newcommand{\edfn}{\end{dfn}} 
\newcommand{\ethm}{\end{thm}} 
\newcommand{\elmma}{\end{lmma}} 
\newcommand{\eppsn}{\end{ppsn}} 
\newcommand{\ecrlre}{\end{crlre}} 
\newcommand{\exmpl}{\end{xmpl}} 
\newcommand{\ermrk}{\end{rmrk}} 
 
\newcommand{\IC}{{C\!\!\!|~}}

\newcommand{\al}{\alpha}

\newcommand{\cla}{{\cal A}} 
\newcommand{\clb}{{\cal B}} 
\newcommand{\clc}{{\cal C}} 
 
\newcommand{\cle}{{\cal E}} 
\newcommand{\clf}{{\cal F}} 
 
\newcommand{\clh}{{\cal H}}

\newcommand{\clk}{{\cal K}} 
\newcommand{\cll}{{\cal L}} 
\newcommand{\clm}{{\cal M}}

\newcommand{\clt}{{\cal T}}

\def\a*{{\cal A}_{h,*}} 
\def\B{{\cal B}(h)} 
\def\B1{{\cal B}_1(h)} 
\def\b{{\cal B}^{s. a. }(h)} 
\def\b1{{\cal B}^{s. a. }_1(h)}

\def\a0{{\cal A}_0}
\def\c0{{\cal C}_0}

\newcommand{\ot}{\otimes}

\newcommand{\raro}{\rightarrow}

\newcommand{\id}{\mbox{id}}

\newcommand{\qed} { \mbox{}\hfill \vspace{1ex}}

\begin{document} 
 \begin{center}
{\large {\bf Towards the Baum-Connes' Analytical Assembly Map for the Actions
of Discrete Quantum Groups}}\\
by\\
{\large Debashish Goswami}\\
and\\
{\large A. O. Kuku;}\\
{\large The Abdus Salam International Centre for Theoretical Physics
(Mathematics Section)}\\ 
{\large Strada Costiera 11, trieste 34014, Italy.}\\
{\large e-mail : goswamid@hotmail.com,}\\
{\large kuku@ictp.trieste.it}\\
Address for correspondence :\\
{\large Prof. A. O. Kuku, Mathematics Section, The Abdus Salam I.C.T.P.}\\
{\large Strada Costiera 11, Trieste 34014, Italy.}\\
{\large Phone : 0039 - 040 - 2240267, Fax :0039 - 040 - 224163}\\
{\large Running head : {\bf Baum-Connes Map for Quantum Groups}}\\
\end{center}
\newpage
\begin{abstract}
Given an action of a discrete quantum group (in the sense of Van Daele,
Kustermans and Effros-Ruan) $\cla$ on a $C^*$-algebra $\clc$, satisfying some
 regularity assumptions resembling the proper $\Gamma$-compact action for a
classical discrete group $\Gamma$ on some space, we are able to construct
canonical maps $\mu^i_r$  ($\mu^i$ respectively)  ($i=0,1$ ) from the
$\cla$-equivariant K-homology groups $KK_i^\cla(\clc,\IC)$ 
   to the K-theory groups $K_i(\hat{\cla}_r)$  ($K_i(\hat{\cla})$
respectively), where $\hat{\cla}_r$ and $\hat{\cla}$ stand for the quantum
analogues of the reduced and full group $C^*$-algebras (c.f. \cite{van2},
\cite{er}). We follow the steps of the construction of the classical
Baum-Connes map (c.f. \cite{bc},\cite{bch}, \cite{val1}, \cite{val2}), although
in the context of quantum group the nontrivial modular property of the
invariant weights (and the related fact that the square of the antipode is not
identity) has to be taken into serious consideration, making it somewhat
tricky to guess and prove the correct definitions of relevant Hilbert module
structures. 
\end{abstract}
\vspace{7mm}
\newpage
  \section{Introduction}
The famous conjecture  made by Paul Baum and Alain Connes has given
birth to one of the most interesting areas of research in both classical and
noncommutative geometry, topology, K-Theory etc. Let us very briefly recall
the main statement of this conjecture (c.f. \cite{bc}, \cite{bch},
and also \cite{val1}, \cite{val2} for a nice and easily accessible account).
Given a locally compact group $G$, and a locally compact Hausdorff space
$X$ equipped with a $G$-action such that $X$ is  proper and $G$-compact (see
for example \cite{val1} and the references therein for various equivalent
formulation of these concepts), there are canonical maps $\mu^r_i :
KK_i^G(C_0(X),\IC) \raro K_i(C^*_r(G))$, and $\mu_i : KK_i^G(C_0(X),\IC) \raro
K_i(C^*(G)),$ for $i=0,1,$ where $C_0(X)$ is the commutative $C^*$-algebra of
continuous complex-valued functions on $X$ vanishing at infinity, $C^*_r(G)$
and $C^*(G)$ are respectively the reduced and free groups $C^*$-algebras, and
$KK^G_.$ denotes the Kasparov's equivariant KK-functor. In particular,
$KK^G(C_0(X),\IC)$ is identified with the $G$-equivariant K-homology of $X$,
and thus is essentially something geometric or topological, whereas the object
 $K_i(.)$ on the right hand side involves the reduced or free group algebras,
which are analytic in some sense. Now, let $\underline{EG}$ be the universal
space for proper actions of $G$. The definition of proper $G$-actions  and
explicit constructions in various cases of interest can be found in
\cite{val1},\cite{val2} and the references therein. The equivariant
$K$-homology of $\underline{EG}$, say $RK^G_i(\underline{EG}),i=0,1,$ can be
defined as the inductive limit of $KK^G_i(C_0(X),\IC),$ over all possible
locally compact, $G$-proper and $G$-compact subsets $X$ of the universal
space $\underline{EG}.$  Since the construction of $KK^G_i$ and $K_i$ commute
with the procedure of taking an inductive limit, it is possible to define
$\mu^r_i,\mu^i$ on the equivariant $K$-homology $RK^G_i(\underline{EG})$, and
the conjecture of Baum-Connes states that $\mu^r_i,i=0,1$ are isomorphisms of
abelian groups. This conjecture admits certain other generalizations, such as
the Baum-Connes conjecture with coefficients (which seems to be false from
some recent result announced by M. Gromov, see \cite{val1} for references),
but we do not want to discuss those here. However, we would like to point out
that the Baum-Connes conjecture has already been verified for many classical
groups, using different methods and ideas from many diverse areas of
mathematics, and has given birth to many new and interesting tools and
techniques in all these areas.  In fact, the truth of this conjecture, if
established, will prove many other famous conjectures in topology, geometry
and $K$-theory. 

Now, in last two decades, the theory of quantum groups has become another
fast-growing branch of mathematics and mathematical physics. Motivated by
examples coming from physics, as well as some fundamental mathematical problems
 (e.g. to develop a good theory of duals for noncommutative topological
groups), many mathematicians including Drinfeld, Jimbo, Woronowicz and others
have formulated and studied the concept of quantum groups, which is a
far-reaching generalization of classical topological groups. On the other
hand, with the pioneering   efforts of Connes (see \cite{Con}), followed by
himself and many other mathematicians, a powerful generalization of classical
differential and Riemannian geometry has emerged under the name of
noncommutative geometry, which has had, since its very beginning, very close
connections with $K$-theory too. Furthermore, 
Baaj and Skandalis (\cite{bs}) have been able to construct an analogue of
equivariant $KK$-theory for the actions of quantum groups, 
 as natural extension of Kasparov's equivariant $KK$-theory. This motivates
one to think of a possibility of generalizing the Baum-Connes construction in
the framework of quantum groups. In the present article, we make an 
attempt towards this generalization.  As we have already mentioned : there are
two steps in  the classical formulation of the Baum-Connes conjecture.  
  First of all, one has to define the maps $\mu^r_i$ for
$G$-compact and $G$-proper actions. Then in the second step, one defines the
universal space for proper action of the group, and then more importantly,
tries to build explicit good models for this universal space to show that it
can be approximated in a suitable sense by its subsets having $G$-proper and
$G$-compact actions, thereby defining the maps $\mu^r_i$ by inductive limit.
What we have been able to achieve in our work here is essentially the first
step, for a class of quantum groups called discrete quantum groups (which are
indeed generalizations of  discrete groups). However, a definition of
proper action of quantum groups has been already proposed in \cite{prop}, and
we hope that it may be possible to achieve the second step starting from this
definition, thereby actually formulating (and then verifying in some cases, if
possible) Baum-Conens conjecture for discrete quantum groups. But we would
like to postpone that task for later work.  

We would also like to mention one thing. We have restricted ourselves within
the framework of discrete quantum groups not only because it is technically
easier to do so, but also because, in fact, the classical Baum-Connes
conjecture is very interesting and nontrivial for discrete groups, and in some
sense most of the difficult cases belong to them. Of course, if our programme
seems to go through satisfactory for discrete quantum groups, we would like to
take up more general locally compact quantum groups in future. It should be
noted that  in the quantum case, discreteness does not imply the unimodularity
of the haar weight, and thus even for discrete quantum groups, one has to be
very careful about the choices of left or right invariant weights as well as
the  appropriate role of the modular operator, as we shall see.

Let us conclude this section with  some useful notational convention. For a
Hilbert space $\clh$, and some pre-$C^*$-algebra $\clb \subseteq \clb(\clh),$
we shall denote the multiplier algebra of the norm-closure of $\clb$ by
$\clm(\clb).$ For two Hilbert spaces $\clh_1,\clh_2$ and some bounded operator
$X \in \clb(\clh_1 \ot \clh_2)=\clb(\clh_1) \ot \clb(\clh_2),$ we denote by
$X_{12}$ the operator $ X \ot 1_{\clh_2}$ on $\clh_1 \ot \clh_2 \ot \clh_2,$
and denote by $X_{13}$ the operator $(1_{\clh_1} \ot \Sigma) (X \ot
1_{\clh_2})(1_{\clh_1} \ot \Sigma)$ on $\clh_1 \ot \clh_2 \ot \clh_2$, where
$\Sigma : \clh_2 \ot \clh_2 \raro \clh_2 \ot \clh_2$ flips the two copies of
$\clh_2.$ For two vectors $\xi,\eta \in \clh_1$ we define a map $T_{\xi \eta}
: \clb(\clh_1 \ot \clh_2) \raro \clb(\clh_2)$ by setting $T_{\xi \eta}(A \ot
B):=<\xi,A\eta>B,$ where $A \in \clb(\clh_1), B \in \clb(\clh_2),$ and extend
this definition to the whole of $\clb(\clh_1 \ot \clh_2)$ in the obvious way. 
It is easy to see that $T_{\xi \eta}(X^*)=(T_{\eta \xi}(X))^*,$ and $T_{\xi
\xi}(X)$ is nonnegative operator if $X$ is. In fact, $T_{\xi \xi}$ is a
completely positive map.  

 For some  Hilbert space $\clh$, we
denote by  $\clb_0(\clh)$  the $C^*$-algebra of compact operators on
$\clh$, and by $\cll(E)$  the $C^*$-algebra of adjointable linear maps on
a Hilbert $\cla$-module $E.$ Furthermore, for a von Neumann algebra $\clb
\subseteq \clb(\clh),$ and some Hilbert space $\clh^\prime,$ we introduce the
following notation : for $\eta \in \clh^\prime, X \in  \clb(\clh^\prime) \ot
\clb \equiv \cll(\clh^\prime \ot \clb),$ $X \eta := X(\eta
\ot 1_\clb) \in \clh^\prime \ot \clb.$ Note that we have
denoted by $\clh^\prime \ot \clb$ the Hilbert von Neumann module obtained from
the algebraic $\clb$-module  $\clh^\prime \ot_{\rm alg} \clb $ by completing this algebraic module
in the strong operator topology inherited from $\clb(\clh,\clh^\prime
\ot \clh),$  where we have identified an element of the form
$(\xi \ot b)$, $\xi \in \clh^\prime, b \in \clb,$ with the
operator which sends a vector $v \in \clh$ to $(\xi \ot bv) \in
\clh^\prime \ot \clh.$  It is easy to see that $\clh^\prime \ot \clb$ is
isomorphic as a Hilbert von Neumann module with $ \{ X \in
\clb(\clh, \clh^\prime  \ot \clh) : X c=(1 \ot c)X, \forall c \in \clb^\prime
\},$ where $\clb^\prime$ denotes the commutant of $\clb$ in $\clb(\clh).$
 Similarly, for a possibly nonunital $C^*$-algebra $\cla$, we can complete the
algebraic pre-Hilbert $\cla$ module $\clh^\prime \ot_{\rm alg} \cla$ in the
locally convex topology coming from the strict topology on $\clm(\cla)$, so
that the completion becomes in a natural way a locally convex Hilbert
$\clm(\cla)$-module, to be denoted by $\clh^\prime \ot \clm(\cla).$ It is
also easy to see that if $X \in \clm(\clb_0(\clh^\prime) \ot \cla) ,$ $\eta
\in \clh^\prime,$ then we have $X\eta \equiv X(\eta \ot 1) \in \clh^\prime
\ot \clm(\cla).$   

If $\clb_1,\clb_2$ are two von Neumann algebras, $\clh^\prime$ is a Hilbert
space,  and $\rho : \clb_1 \raro \clb_2$ is a normal $\ast$-homomorphism, then
it is easy to show that $(id \ot \rho ) : \clh^\prime \ot_{\rm alg} \clb_1
\raro \clh^\prime \ot_{\rm alg} \clb_2$ admits a unique extension (to be
denoted again by $(id \ot \rho)$) from the Hilbert von Neumann module
$\clh^\prime \ot \clb_1$ to the Hilbert von Neumann module $\clh^\prime \ot
\clb_2.$ Furthermore, one has that $(id \ot \rho)(X \eta)=(id \ot \rho)(X)
\eta$ for $X \in \clb(\clh^\prime) \ot \clb,$ $\eta \in \clh^\prime.$ By very
similar arguments one can also prove that if $\cla_1,\cla_2$ are two
$C^*$-algebras, and $\pi : \cla_1 \raro \cla_2$ is a nondegenerate
$\ast$-homomorphism (hence extends uniquely as a unital strictly continuous
$\ast$-homomorphism from $\clm(\cla_1)$ to $\clm(\cla_2)$), then $(id \ot \pi)
: \clh^\prime \ot_{\rm alg} \cla_1 \raro \clh^\prime \ot_{\rm alg} \cla_2$
admits a unique  extension (to be denoted by the same notation) from
$\clh^\prime \ot \clm(\cla_1)$ to $\clh^\prime \ot \clm(\cla_2),$ which is
continuous in the locally convex topologies coming from the respective strict
topologies. We also have that $(id \ot \pi)(X \eta)=(id \ot \pi)(X) \eta,$ for
$X \in \clm(\clb_0(\clh^\prime) \ot \cla), \eta \in \clh^\prime.$       

\section{Preliminaries on discrete quantum groups} 
We briefly discuss the theory of {\it discrete quantum groups} as developed in
\cite{van}, \cite{er},\cite{van3}, \cite{kus} and other relevant references
to be found there. Let us fix an index set $I$ (possibly uncountable), and let
$\a0:=\oplus_{\al \in I} \cla_\al$ be the {\it algebraic} direct sum of
$\cla_\al$'s, where for each $\al,$  $\cla_\al=M_{n_\al}$ is the finite
dimensional $C^*$-algebra of $n_\al \times n_\al$ matrices with complex
entries, and $n_\al$ is some positive integer. Let us denote by $\clm\equiv
\clm(\a0)$ the unital $C^*$-algebra consisting of all collections
$(a_\al)_{\al \in I}$ with $a_\al \in \cla_\al$ for each $\al$, and $\sup_\al
\| a_\al \| < \infty.$ The algebra operations are taken to be the obvious
ones; i.e. $(a_\al)+(b_\al):=(a_\al+b_\al)$, $(a_\al).(b_\al):=(a_\al b_\al)$
and $(a_\al)^*:=(a_\al^*).$   Similarly, denote by $\clm(\a0 \ot \a0)$ the
$C^*$-algebra consisting of all collections of the form $(a_\al \ot b_\beta)$
where $\al, \beta $ varies over $I$, and equip $\clm(\a0 \ot \a0)$ with the
obvious $C^*$-algebra structure. Let us now assume that there is a unital
$C^*$-homomorphism $\Delta : \clm(\a0) \raro \clm(\a0 \ot \a0)$ which
satisfies the following :\\  (i) For $a , b \in \a0,$ we have $$ T_1(a \ot b)
:=\Delta(a)(1 \ot b) \in \a0 \ot_{\rm alg} \a0,$$  and
$$ T_2(a \ot b) :=(a \ot 1) \Delta(b) \in \a0 \ot_{\rm alg} \a0;$$
(ii) $T_1, T_2 : \a0 \ot_{\rm alg} \a0 \raro \a0 \ot_{\rm alg} \a0$ are
bijections;\\
(iii) $\Delta$ satisfies the coassociativity in the sense that 
$$ (a \ot 1 \ot 1)(\Delta \ot id)(\Delta(b)(1 \ot c))=(id \ot \Delta)((a \ot
1)\Delta(b))(1 \ot 1 \ot c),$$
 for $a,b, c \in \a0.$ 
As explained in the relevant references mentioned above, $(\Delta \ot id),$
 $(id \ot \Delta)$ 
admit  extensions as  $C^*$-homomorphisms from $\clm(\a0 \ot \a0 )$ to
$\clm(\a0 \ot \a0 \ot \a0)$ (we denote these extensions by the same
notation) and the condition (iii) translates into $(\Delta \ot id) \Delta=(id
\ot \Delta) \Delta.$ 

The above conditions essentially constitute the definition of a discrete
quantum group (for details  see \cite{van}, \cite{kus}  and \cite{er}).  Let us
recall from \cite{van} and \cite{kus} some of the important properties of our
discrete quantum group $\a0.$ It is remarkable that it is possible to deduce
from (i) to (iii) the existence of a canonical antipode $S :\a0 \raro \a0$
satisfying $S(S(a)^*)^*=a$ and other usual properties of the antipode of a
Hopf algebra. Furthermore, there exists a counit $\epsilon :\a0 \raro \IC.$
For details of the constructions of these maps and their properties we refer
to \cite{van}.  

\vspace{5mm}
We shall call an  arbitrary collection $(a_\al)_{\al \in I},$ with
$a_\al \in \cla_\al \forall \al$, the ``algebraic multiplier" of $\a0.$ 
 The set of all algebraic multipliers of $\a0$, denoted  by $\clm_{\rm alg}(\a0)$, is
obviously a $\ast$-algebra, with pointwise multiplication and adjoint, i.e.
for $U =(u_\al),V=(v_\al) \in \clm_{\rm alg}(\a0),$ $UV:=(u_\al v_\al)_\al$, and
$U^*:=(u_\al^*).$  Clearly, any element of $\a0$ can be viewed as an element
of $\clm_{\rm alg}(\a0)$, by thinking of $a \in \a0$ as $(a_\al)$, where $a_\al$ is
the component of $a$ in $\cla_\al.$ It is easy to see that $Ua,aU \in \a0$
for $a \in \a0, U \in \clm_{\rm alg}(\a0).$ We can give a similar definition of
algebraic multiplier of $\a0 \ot \a0$, which will be any collection of the
form $M \equiv (m_{\al \beta})_{\al,\beta \in I}$, with $m_{\al \beta} \in
\cla_\al \ot \cla_\beta.$ In fact,
 since by \cite{van} there is a bijection of the index set $I$, say $\al
\mapsto \al^\prime$, such that $S(e_\al)=e_{\al^\prime},
S(e_{\al^\prime})=e_\al,$ we can define $S(X)$ for $X =(x_\al)_I \in
\clm_{\rm alg}(\a0)$ by $S(X):=X^\prime=(x^\prime_\al)_I$ where
$x^\prime_\al=S(x_{\al^\prime}).$ Similarly, to define $\Delta(X)$ for
$X=(x_\al) \in \clm_{\rm alg}(\a0),$ we note that (c.f. \cite{van}) for fixed
$\al,\beta \in I$, there is a finite number of $\gamma \in I$ such that
$\Delta(e_\gamma)(e_\al \ot e_\beta)$ is nonzero. Thus, $\Delta(X)$ can be
defined as the element $Y \in \clm_{\rm alg}(\a0 \ot \a0)$ such that $Y=(y_{\al
\beta}),$ where $y_{\al \beta}= \sum_\gamma \Delta(x_\gamma)
\Delta(e_\gamma)(e_\al \ot e_\beta).$    For algebraic multipliers $A,B$ of $
\a0$ and $L$ of $\a0 \ot \a0,$ it is clear that $\Delta(A)=L$ if and only if 
$\Delta(Aa)=L\Delta(a)$ $\forall a \in \a0,$ and  $S(A)=B$ if and only if 
$S(Aa)=S(a)B$ for $a \in \a0.$ 

\vspace{5mm}
Let $\clk$ be the smallest Hilbert space containing the algebraic direct sum
$\oplus_{\al \in I} \clk_\al \equiv \oplus_\al \IC^{n_\al},$ i.e. $\clk=\{
(f_\al)_{\al \in I} : f_\al \in \clk_\al =\IC^{n_\al}, \sum_\al \| f_\al \|^2 <
\infty \},$ where the possibly uncountable sum $\sum_\al$ means the limit over
the net consisting of all  possible sums over finite subsets of $I$. Let us
consider the canonical imbedding of $\a0$ in $\clb(\clk)$, with $\cla_\al$
acting on $\IC^{n_\al}$. Let $\cla$ be the completion of $\a0$ under the
norm-topology inherited from $\clb(\clk).$ Let us fix some  matrix units 
$e^\al_{ij}, i,j =1,...,n_\al$ for $\cla_\al=M_{n_\al}$, w.r.t. some fixed
orthonormal basis $e^\al_i, i=1,...,n_\al,$ of $\IC^{n_\al},$ and thus $\cla$
is the $C^*$-algebra generated by $e^\al_{ij}$'s.   It is also clear that any
element of $\clm_{\rm alg}(\a0)$ can be viewed as a possibly unbounded operator on
$\clk$, with the domain containing the algebraic direct sum of $\clk_\al$'s.
Similarly, elements of $\clm_{\rm alg}(\a0 \ot \a0)$ can be thought of as possibly
unbounded operators on $\clk \ot \clk$ with suitable domain. 

\vspace{5mm} 
Let us denote by $\a0^\prime$ the set of all linear functionals on $\a0$
having ``finite support", i.e. they vanish on 
$\cla_\al$'s for all but finite many $\al \in I.$ It is clear that any $f \in
\a0^\prime$ can be identified as a functional on $\clm_{\rm alg}(\a0)$, by defining
$f((a_\al)_I):=\sum_{\al \in I} f(a_\al) \equiv \sum_{I_0} f(a_\al),$ where
$I_0$ is the finite set of $\al$'s such that for $\al$'s not belonging to
$I_0$, $f|_{\cla_\al}=0.$ With this identification, $f(1)$ makes sense for any
$f \in \clm_{\rm alg}(\a0)$. Let us denote by $e_\al$ the identity of
$\cla_\al=M_{n_\al}$, which is a minimal central projection in $\a0.$ For any
subset $I_1$ of $I$ we denote by $e_{I_1}$ the direct sum of $e_\al$'s for
$\al \in I_1.$ It is clear that a functional $f$ on $\a0$ is in $\a0^\prime$
if and only if there is some finite $I_1$ such that $f(a)=f(e_{I_1}a)$ for all
$a \in \a0.$

\vspace{5mm}
We say that a linear functional $\phi$ (not necessarily with finite
support) on $\a0$ is left invariant if we have $(id \ot
\phi)((b \ot 1)\Delta(a))=b\phi(a)$ for all $a,b \in \a0,$ or equivalently,  
$\phi((\omega \ot id)(\Delta(a)))=\omega(1) \phi(a)$ for all $a \in \a0,$
$\omega \in \a0^\prime.$ Similarly, a linear functional $\psi$ on $\a0$ is
called right invariant if $(\psi \ot id)((1 \ot b)\Delta(a))=\psi(a)b$ for all
$a,b \in \a0.$ Let us now recall some of the main results regarding left and
right invariant functionals as proved in \cite{van}. It is shown in \cite{van}
that up to constant multiples, there is a unique left invariant functional, and
  same thing is true for right invariant functionals, although in general
(unless $S^2=id$) left and right invariant functionals are not the same.
Moreover, for each $\al \in I,$ there is a positive invertible element $K_\al
\in \cla_\al$ such that the positive functional $\phi$ defined by $$
\phi(x)=Tr_\al(K_\al^{-1}x)$$ for $x \in \cla_\al$, (where $Tr_\al$ is the
trace on the  algebra $\cla_\al$ of $n_\al \times n_\al$ matrices) and extended
on $\a0$ by linearity, is left invariant. We get a right invariant positive
functional $\psi$ by replacing $K_\al^{-1}$ by $c_\al K_\al$ for some positive
constant $c_\al$, i.e. $$\psi(x) :=c_\al Tr_\al(K_\al x),$$
 for $x \in \cla_\al.$ Furthermore, $S^2(a)=K_\al^{-1} a K_\al$ for $a \in
\cla_\al.$ If we define a possibly unbounded positive invertible operator $K$
on $\clk$ by setting $K|_{\clk_\al}=K_\al$ for each $\al,$ then it is easy to
see that $\psi(a)=c_\al \phi(K^2 a)=c_\al \phi(aK^2)$ for $a \in \cla_\al.$
Now, observe that $a \mapsto \phi(S(a))$ is right invariant, hence there is
some constant $c$ such that $\psi=c \phi \circ S.$     From the results of
\cite{kus} it follows that there is a ``modular operator" $\delta$, which can
be thought of as a collection $(\delta_\al)_{\al \in I}$ such that $\delta_\al
\in \cla_\al$ for each $\al$ (i.e. $\delta \in \clm_{\rm alg}(\a0)$), and we also have
that  $\Delta(\delta)=\delta \ot \delta,$ 
$S(\delta)=\delta^{-1},$ $S(\delta^{-1})=\delta,$ in the sense 
 described earlier; and furthermore,  $\phi(S(a))= \phi(a \delta_\al)$ for all
$a \in \cla_\al.$ Thus, for $a \in \cla_\al,$ $ c_\al \phi(aK^2_\al)=\psi(a)=c
\phi(S(a))=c \phi(a \delta_\al).$ Since $c$ is clearly  nonzero, we conclude
that $\delta_\al=c^{-1}c_\al K^2_\al$ for each $\al.$ Let us now argue that
$c$ is positive, which will show the positivity of $\delta_\al$. Since $\psi$
is by construction a positive functional, we need to prove that $\phi \circ S$
is positive too. However, for $a \in \a0,$ 
$\phi(S(a^*a))=\phi(S(a)S(a^*))=\phi(S(a)S^2(S(a)^*))=\phi(S(a)K^{-1}S(a)^*K)$. 
From the definition of   $\phi$ in terms of trace on each finite dimensional
component, it is clear that
$\phi(S(a)K^{-1}S(a)^*K)=\phi(K^{\frac{1}{2}}S(a)K^{-1}
S(a)^*K^{\frac{1}{2}}) \geq 0.$ So, $c$ is positive, and hence so is the
oparator $\delta_\al$ for each $\al.$ Let
$\theta_\al:=\delta_\al^{\frac{1}{2}}$ for each $\al,$ and let $\theta$ be the
unbounded positive operator on $\clk$ defined by
$\theta|_{\clk_\al}=\theta_\al.$ 

\vspace{5mm}
Let us fix some $\al$ now.   From \cite{van}, note that there is some
index $\beta$ such that $S(\cla_\al)=\cla_\beta$, and in particular
$S(e_\al)=e_\beta.$ Since $\delta_\al$  is a finite dimensional positive
invertible matrix, all  its eigenvalues are strictly positive. Similar thing
is true for $\delta_\beta,\delta_\al^{-1},\delta_\beta^{-1}$ too.  Thus, we can
choose a holomorphic function $g$ defined on an open set of the complex plane
containing the union of the spectrum of the matrices
$\delta_\al,\delta_\beta,\delta_\al^{-1},\delta_\beta^{-1}$  such that $g(\delta_\al)=\theta_\al,$
$g(\delta_\al^{-1})=\theta_\al^{-1}.$   As the restriction of $S$ on
$\cla_\al$, say $S_\al$, is a linear map on a finite dimensional space, it is
norm-continuous, and furthermore, $S(x^n)=S(x)^n$ for any positive integer
$n$, $x \in \cla_\al,$ from which it is easy to see that
$S(\theta_\al)=S(g(\delta_\al))=g(S(\delta_\al))=g(\delta_\beta^{-1})=\theta_\beta^{-1}=S(e_\al)\theta^{-1}.$  
 Similarly, $S(\theta_\al^{-1})=\theta_\beta.$ Since this is true for any
$\al$, we conclude that $S(\theta)=\theta^{-1}$ and $S(\theta^{-1})=\theta.$ 
 By a very similar argument we can prove that $\Delta(\theta)=\theta \ot
\theta.$ Furthermore, from our discussion it is also clear that
$S^2(a)=\theta^{-1} a \theta$ for $a \in \a0.$ Let us summarize these facts
here :\\ 
(a) There exists a positive (possibly unbounded) invertible operator
$\theta$ on $\clk$, with its domain containing all $\clk_\al$'s, with
$\theta_\al=\theta|_{\clk_\al} \in \cla_\al$ satisfying
$\Delta(\theta)=(\theta \ot \theta),$ $S(\theta)=\theta^{-1},$ and
$S(\theta^{-1})=\theta.$\\ 
(b) $S^2(a)=\theta^{-1}a\theta$ for all $a \in
\a0.$\\ 
(c) We can choose a positive faithful left invariant functional (to be
referred to as left haar measure later on) $\phi$ and a positive faithful
right invariant functional (to be referred to as right haar measure) $\psi$
such that $\psi(a)=\phi(a \theta^2)=\phi(\theta^2 a)$ for $a \in \a0.$\\ (d)
$\phi(S^2(a))=\phi(a), \psi(S^2(a))=\psi(a)$ for all $a \in \a0$, where
$\phi,\psi$ as in (c).\\ Note that we may have to multiply the left and right
invariant functionals $\phi$ and $\psi$ we constructed earlier by some
positive constant in order to  make them satisfy the property (c) above. 

\vspace{5mm}
 We say that a unitary element in $\clm(\cll(\clh \ot \cla)) \equiv
\clm(\clb_0(\clh) \ot \cla)$ is a unitary representation of the discrete
quantum group $\cla$ if $(id \ot \Delta)(U)=U_{12} U_{13},$ and $(id \ot
S)(U)=U^*.$ Note that the second equality has to be understood in the sense
 of the definition of $S$ on the algebraic multiplier, i.e. $(id \ot
S)(U(1 \ot a))=(1 \ot S(a))U^*$ for all $a \in \a0.$ Let us also make the
following useful observation : for $X \in \clm(\clb_0(\clh) \ot \cla),$ and
$\xi, \eta \in \clh,$ we have that $T_{\xi,\eta}(X) \in \clm(\cla).$ 

\vspace{5mm}
Let us now extend the definition of $\phi$ and $\psi$ on a larger set than 
$\a0$ as follows. For a nonnegative element $a \in \clm(\cla) \subseteq
\clb(\clk),$ we define $\phi(a)$ as the limit of $\phi_J(a)$, whenever this
limit exists as a finite number, and  where $J$ is any finite subset of $I$,
$\phi_J(.):=\phi(e_J.)=\phi(.e_J),$ and the limit is taken over the net of
finite subsets of $I$ partially ordered by inclusion. Similarly, we set
$\psi(a)=\lim_J \psi(e_Ja)$ whenever the limit exists as a finite number.
Since  a general element $a \in \clm(\cla)$ can be    canonically written as a
linear combination of four nonnegative elements, and extend the definition of
$\phi$ on $\clm(\cla)$ by linearity.    For any nonnegative $X \in
\clm(\clb_0(\clh) \ot \cla)$ (where $\clh$ is some Hilbert space), we define
$(id \ot \phi)(X)$ as the  limit in the weak-operator topology  (if it exists
as a bounded operator) of the net $(id \ot \phi_J)(X)$ over finite subsets $J
\subseteq I,$ and extend this definition for a general $X \in
\clm(\clb_0(\clh) \ot \cla)$ in the usual way.    Similar definition will be
given for $(id \ot \psi)$.  
\blmma
 If we choose $\clh=\clk$ in the
above, and take any $a \in \clm(\cla)$ such that $\phi(a)$ is finite, then
$(id \ot \phi)(\Delta(a))=\phi(a)1_{\clm(\cla)}.$ 
\elmma
{\it Proof :-}\\
For any nonnegative $X \in \clm(\cla \ot \cla)$, and any positive operator $P
\in \a0 \ot \a0$, with $0 \leq P \leq 1$, such that $P$ and $X$ commute, it is
easy to see that $(id \ot \phi_J)(PX) \leq (id \ot \phi_J)(X)$, for any finite
subset $J$ of $I.$  By choosing large enough $J$ one can ensure that $P
\leq  (1 \ot e_J)$, so that $(id \ot \phi_J)(PX)=(id \ot \phi)(PX)$ ($PX$ is
in $\a0 \ot \a0$, so $(id \ot \phi)(PX)$ makes sense). So, for any vecror $\xi
\in \clk,$ $\sup_P <\xi, (id \ot \phi)(PX) \xi> \leq \sup_J <\xi, (id \ot
\phi_J)(X) \xi>,$ where the supremum in the left hand side is taken over all
positive $P \in \a0 \ot \a0$ with $P \leq 1,$ and commutes with $X$.  On the
other hand, for fixed finite subsets $J,K$, $(e_K \ot e_J)$ is one such $P$,
and thus $  \sup_P <\xi, (id \ot \phi)(PX) \xi> \geq <e_K \xi, (id \ot
\phi_J)(X) e_K \xi>,$ and taking limit over $K$, we conclude that 
$\sup_P <\xi, (id \ot \phi)(PX) \xi> \geq \sup_J <\xi, (id \ot
\phi_J)(X) \xi>,$
 which proves that they are equal, and hence $(id \ot \phi)(X)$ exists as a
bounded operator if and only if the weak-operator limit of $(id \ot
\phi)(P_\nu X)$ exists as a bounded operator, over any net $P_\nu$ of
nonnegative operators in $\a0 \ot \a0$, commuting with $X$, and such that
$P_\nu \uparrow 1.$   Using this fact, we see that for nonnegative $a \in
\clm(\cla)$, $(id \ot \phi)(\Delta(a))=\lim_{J,K} (id \ot \phi)((e_K \ot
1)\Delta(e_J)\Delta(a))=\lim_{J,K} (id \ot \phi)((e_K \ot 1)\Delta(ae_J)),$
where $J,K$ are varied over all finite subsets of $I$, and we have used the
fact that $e_K$'s are central projections, and thus $\Delta(e_J)$ commutes
with $\Delta(a)$. From the above expression, by using the left invariance of
$\phi$ on $\a0$, and then taking limits, the desired result follows.
 \qed 

We remark that an analogous fact is true for $\psi.$ 

\vspace{5mm}
We shall now define a $\ast$-algebra structure on $\a0^\prime$, and then
identify $\a0$ with suitable elements of $\a0^\prime,$ thereby equipping $\a0$
with this new $\ast$-algebra structure, and then consider suitable
$C^*$-completions.  This will give rise to the analogues of the full and
reduced group $C^*$-algebra in the framework of discrete quantum groups.
Following \cite{kus} and others, we define $f \ast g $ for $f,g \in \a0$ by
$(f \ast g)(a):=(f \ot g)(\Delta(a)), a \in \a0.$ Note that since $f,g$ have
finite supports, there is some finite subset $J$ of $I$ such that $(f \ot
g)(\Delta(a))=(f \ot g)((e_J \ot e_J)\Delta(a)),$ and since $(e_J \ot
e_J)\Delta(a) \in \a0 \ot_{\rm alg} \a0,$ $f \ast g$ is well defined. We also
define an adjoint by $f^*(a):=\bar{f}(S(a)^*), a \in \a0.$ We now define for
each $a \in \a0$, an element $\psi_a \in \a0^\prime$ by $\psi_a(b):=\psi(ab).$ 
 It is easy to verify the following by using standard formulae involving
$\Delta$ and $S.$
\bppsn
For $a,b \in \a0,$ $\psi_a \ast \psi_b=\psi_{a \ast b},$ where $a \ast b :=(id
\ot \psi)((1 \ot b)((id \ot S^{-1})(\Delta(a))))=(\phi
\ot id)((a \ot 1)((S \ot id)(\Delta(b)))).$ Furthermore,
$\psi_a^*=\psi_{a^\sharp},$ where $a^\sharp:=\theta^{-2}S^{-1}(a^*).$ 
\eppsn   

We denote by $\hat{\a0}$ the set $\a0$ equipped with the $\ast$-algebra
structure given by $(a,b) \mapsto a \ast b, a \mapsto a^\sharp$ described by
the above proposition. There are two different natural ways of making
$\hat{\a0}$ into a $C^*$-algebra, and thus we obtain the so-called reduced
$C^*$-algebra $\hat{\cla}_r$ and the free or full $C^*$-algebra $\hat{\cla}$.
This is done in a similar way as in the classical case : one can realize
elements of $\hat{\a0}$ as bounded linear operators on the Hilbert space
$L^2(\phi)$ (the GNS-space associated with the positive linear functional
$\phi$, see \cite{van2} and \cite{er} for details) and complete $\hat{\a0}$ in
the norm inherited from the operator-norm of $\clb(L^2(\phi))$ to get
$\hat{\cla}_r$. The definition of $\hat{\cla}$ is slightly more complicated
and involves the realization of $\hat{\a0}$ as elements of the
 Banach $\ast$-algebra $L^1(\phi)$ (see \cite{er} and other relevant
references) and then taking the associated universal $C^*$-completion. 
However, it is not important for us how the explicit constructions of these two
$C^*$-algebras are done; we refer to \cite{van2}, \cite{er} for that; all we
need is that  $\hat{\a0}$ is dense in both of them in the respective
norm-topologies.    It should also be mentioned that exactly as in the
classical case, there is a canonical surjective $C^*$-homomorphism from
$\hat{\cla}$ to $\hat{\cla}_r.$

 \newsection{Construction of the analytic assembly map}
In this section, we shall show how one can construct an analogue of the
Baum-Connes analytic assembly map for the action of the discrete quantum group
$\a0$ on some $C^*$-algebra, under some additional assumptions on the action,
which may be called ``properness and $\cla$-compactness", since these
assumptions are actually weaker than having a proper and $G$-compact action in
 the classical situation of an action by a group $G$. Our construction is 
analogous   to that described in, for example, \cite{val1},\cite{val2}, for the
  discrete group. We essentially translate that into our noncommutative
framework step by step, and verify that it really goes through. However, in
case $S^2$ is not identity, it is somewhat tricky to give the correct
definition of $\hat{\a0}$-valued inner product, and prove the required
properties, as one has to suitably incorporate the modular operator $\delta$. 
  
\vspace{5mm}
Let $\clc $ be a $C^*$-algebra (possibly nonunital).
Assume furthermore that there is an action of the quantum group  $\cla$ on it,
given by $\Delta_\clc :\clc \raro \clm(\clc \ot \cla)$, which is coassociative
$C^*$-homomorphism, and assume also that there is a dense $\ast$-subalgebra
$\clc_0$ of $\clc$ such that the following conditions are satisfied :\\
{\bf A1} $\Delta_\clc(c)(c^\prime \ot 1) \in \clc_0 \ot_{\rm alg} \a0$ for all
$c,c^\prime \in \clc_0;$ \\ 
{\bf A2} $\Delta_\clc (c)(1 \ot a) \in \clc_0 \ot_{\rm alg} \a0$ for all $c
\in \clc_0,$ $a \in \a0;$\\
{\bf A3} There is a positive element $h \in \clc_0$ such that $$(id \ot
\phi)(\Delta_\clc(h^2))=1,$$ or equivalently $(id \ot \phi)(\Delta_\clc(h^2)(c
\ot 1))=c, \forall c \in \clc_0 .$

\brmrk
In the classical situation, when $\cla$ is $C_0(G)$ for some discrete group,
and $\a0=C_c(G)$, $\clc=C_0(X)$ for some locally compact Hausdorff space $X$
equipped with an $G$-action such that $X$ is $G$-compact and $G$-proper,
  one can take $\clc_0=C_c(X)$, and it is easy to verify that with this
 choice of $\clc_0$, the conditions {\bf A1}, {\bf A2}, {\bf A3} are satisfied.
It is rather stratightforward to see {\bf A1} and {\bf A2}. The construction
of a positive function $h$ satisfying {\bf A3} can be found in \cite{val1}.
 Thus, the above conditions are in some sense noncommutative generalization of
$G$-proper and $G$-compact actions. 
\ermrk

\brmrk
\label{subgp}
Let us note that the above assumptions are indeed satisfied in a typical
situation, namely for nice ``quantum quotient spaces" corresponding to
``quantum subgroups" of the discrete quantum group $\cla.$ Indeed, from
\cite{van}, the existence of an element $h$ as in {\bf A3} above follows, if
we take $\clc=\cla,$ and $\clc_0=\a0.$ The same thing will trivially hold if
$\clc$ is taken to be a direct sum of finitely many copies of $\cla$, with the
 natural action of $\cla.$ More importantly, there is a natural generalization
of the notion of subgoups and quotient spaces for quantum groups, which is
 by now more or less well-konwn and standard in this theory (see, for example,
\cite{pod} for these concepts in the context of compact quantum groups, and
note that for more general quantum groups they can be easily extended). It is
not difficult to see, by using the fact that our assumptions {\bf A1,A2,A3} are
valid for $\clc=\cla,\clc_0=\a0$, that the same thing will be true if we take
$\clc$ to be the quotient by some compact (i.e. finite dimensional in this
case) quantum subgroup of $\cla$. This is of particular interest in view of
the fact that one of the models for the universal  space for
proper actions of a classical discrete second countable group involves  some
kind of ``infinite join" of some set constructed out of disjoint union over
all possible quotient spaces by finite subgroups of the group. Thus, if a
similar construction can be done   in the noncommutative framework  starting
from the definition of proper actions as proposed in \cite{prop}, then it
seems very likely that using the techniques of the present article an analytic
assembly map can be defined on the ``quantum unversal  space" for ``quantum
proper action", and hence a precise formulation of the Baum-Connes conjecture
for  discrete quantum groups will turn into a reality. 
    \ermrk

\vspace{5mm}
Now, our aim is to construct maps $\mu_i : KK_i^\cla(\clc,\IC) 
 \raro KK_i(\IC,\hat{\cla}) \equiv K_i(\hat{\cla}),$ and $\mu_i^r :
KK_i^\cla(\clc,\IC)   \raro KK_i(\IC,\hat{\cla}_r) \equiv
K_i(\hat{\cla}_r),$ for $i=0,1$, i.e. even and odd cases. For simplicity let
us do it for $i=1$ only, the other case can be taken care of by obvious
modifications. We have chosen the convention of \cite{val1} to treat
separately odd and even cases, instead of treating both of them on the same
footing as in the original work of Kasparov or in \cite{kk}. This is merely
a matter of notational simplicity. For the definition and properties of
equivariant KK groups $KK^\cla_.(.,.)$, we refer to  the paper by Baaj and
Skandalis (\cite{bs}) (with the easy modifications of their definitions to
treat odd and even cases separately). 

Let $(U,\pi,,F)$ be a
cycle (following \cite{val1}) in $KK_1^\cla(\clc,\IC),$ i.e.\\ 
  (i) $U \in 
\cll(\clh \ot \cla)) \cong \clm(\clb_0(\clh) \ot \cla)$ is a unitary
representation of $\cla$, where  $\clh$ is a separable Hilbert space, i.e. $U$
is unitary and $(id \ot \Delta)(U)=U_{12}U_{13},$ $(id \ot S)(U)=U^*;$ \\    
(ii) $\pi : \clc \raro \clb(\clh)$ is a nondegenerate $\ast$-homomorphism such
that $(\pi \ot id)(\Delta_\clc(a))=U(\pi(a) \ot 1)U^*, \forall a \in \clc$;\\ 
(iii) $F \in \clb(\clh)$ is  self-adjoint,  $[F,\pi(c)],\pi(c)(F^2-1) \in
\clb_0(\clh) \forall c \in \clc,$  and $(F \ot 1)-U(F \ot 1)U^* \in
\clb_0(\clh) \ot \cla$.

We say that a cycle $(U,\pi,F)$ is equivariant (or $F$ is equivariant) if $U(F
\ot 1)U^*=F \ot 1$. We say that $F$ is properly supported if for any $c \in
\clc_0,$ there are {\bf finitely many} $c_1,...c_k, b_1,...,b_k \in \clc_0$
and $A_1,...A_k \in \clb(\clh)$ (all depending on $c$) such that $F
\pi(c)=\sum_i \pi(c_i)A_i \pi(b_i).$    

Before we proceed further, let us make the following convention : we
canonically embed $\cla$ in the set of bounded operators on  $ \clk,$ as
described before, and for any element $A \in \clb(\clh),$ we shall denote by
$\tilde{A}$ the element $A \ot 1_{\clk}$ in $\clb(\clh \ot \clk)$.  

\bthm
\label{proper} 
Given a cycle $(U, \pi, F)$, 
 we can find  a homotopy-equivalent  cycle $(U,\pi,
F^\prime)$ such that $(U,\pi, F^\prime)$ is equivariant and $F^\prime$ is
properly supported.
\ethm
{\it Proof :-}\\
 Since $\pi$ is nondegenerate, we can
choose a net $e_\nu$ of elements from $\clc_0$ such that $\pi(e_\nu)$
converges to the identity of $\clb(\clh)$ in the strict topology, i.e. in the
strong $\ast$-topology. Now, let $X_\nu:=\tilde{\pi(e_\nu)^*}U
\tilde{\pi(h)}\tilde{F}
\tilde{\pi(h)}U^*\tilde{\pi(e_\nu)}=\tilde{\pi(e_\nu)^*} (\pi \ot
id)(\Delta_\clc(h))U\tilde{F}U^* (\pi \ot
id)(\Delta_\clc(h))\tilde{\pi(e_\nu)}.$ Since by our assumption
$\tilde{e_\nu^*}\Delta_\clc(h) \in \clc_0 \ot_{\rm alg} \a0,$ and similar
thing is true for $\Delta_\clc(h)\tilde{e_\nu}$, it is easy to see that
$X_\nu$ is of the form $X_\nu=\sum_j (\pi(c_j) \ot
a_j)(U\tilde{F}U^*)(\pi(c^\prime_j) \ot a^\prime_j),$ for some finitely many
$c_j, c_j^\prime \in \clc_0$ and $a_j, a^\prime_j \in \a0.$ Choosing a
suitably large enough finite subset $I_1$ of $I$, we can assume that all the
$a_j, a^\prime_j$'s are in the support of $e_{I_1}$, and hence it is easy to
see that $X_\nu \in \clb(\clh) \ot_{\rm alg} (e_{I_1}\a0 e_{I_1}),$ so $(id
\ot \phi)(X_\nu)$ is finite. Similarly, $(id \ot
\phi)(\tilde{\pi(e_\nu)^*}U\tilde{\pi(h^2)}U^*\tilde{\pi(e_\nu)})$ is finite,
and by assumption  {\bf A3}, is equal to $(\pi(e_\nu^* e_\nu) \ot 1).$ Now,
from the operator inequality $-\|F\| 1 \leq F \leq \| F\| 1,$ we get the
operator inequality
$$ -\tilde{\pi(e_\nu)^*}U\tilde{\pi(h^2)}U^*\tilde{\pi(e_\nu)}\|F\| \leq 
 X_\nu \leq \tilde{\pi(e_\nu)^*}U\tilde{\pi(h^2)}U^*\tilde{\pi(e_\nu)}\|F\|;$$
 from which it follows after applying $(id \ot \phi)$ that
$$ -\pi(e_\nu^*e_\nu) \|F \| \leq (id \ot \phi) X_\nu \leq \pi(e_\nu^*e_\nu)
\|F\|.$$
 Since $\pi(e_\nu^*e_\nu) \raro 1_{\clb(\clh)}$ in the strong operator
topology, one can easily prove by the arguments similar to those in
\cite{val1} that $(\id \ot \phi)(X_\nu)$ converges in the strong operator
topology of $\clb(\clh)$, and let us denote this limit by $F^\prime.$  
It is also easy to see that in fact $F^\prime=(id \ot \phi)(U (\pi(h)F\pi(h)
\ot 1)U^*),$ where we have used the extended definition of $(id \ot \phi)$ on
$\clm(\clb_0(\clh) \ot \cla)$ as discussed in the previous section. 

Fix some $c \in \clc_0.$ Clearly we have $F^\prime \pi(c)=(id \ot \phi)(U
\tilde{\pi(h)}\tilde{F} \tilde{\pi(h)}U^* \tilde{\pi(c)}).$ Now, note that 
 $U \tilde{\pi(h)}\tilde{F} \tilde{\pi(h)}U^* \tilde{\pi(c)}=(\pi \ot
id)(\Delta_\clc(h))U\tilde{F}U^* (\pi \ot id)((\Delta_\clc(h)(c \ot 1)).$
Since $\Delta_\clc(h)(c \ot 1) \in \clc_0 \ot_{\rm alg } \a0,$ we can write it
as a finite sum of the form $\sum_{ij,\al} x^\al_{ij} \ot e^\al_{ij},$
with $x^\al_{ij} \in \clc_0,$ and where $e^\al_{ij}$ 's are the matrix units of
$\cla_\al$, as described in the previous section, and $\al$ in the above sum
varies over some finite set $T$, say, with $i,j=1,...,n_\al.$ Thus, $  U
\tilde{\pi(h)}\tilde{F} \tilde{\pi(h)}U^* \tilde{\pi(c)} =\sum_{\al,i,j} (\pi
\ot id)(\Delta_\clc(h))(1 \ot e^\al_{ij})(Fx^\al_{ij} \ot 1).$ Since for each
$\al,i,j,$ $\Delta_\clc(h)(1 \ot e^\al_{ij}) \in \clc_0 \ot_{\rm alg} \a0,$
we can write $\Delta_\clc(h)(1 \ot e^\al_{ij}) $ as a finite sum of the form
$ \sum x_p \ot a_p$ with $x_p \in \clc_0, a_p \in \a0,$ and hence 
$  U \tilde{\pi(h)}\tilde{F} \tilde{\pi(h)}U^* \tilde{\pi(c)}$ is clearly a
finite sum of the form $\sum_k \pi(c_k)A_k\pi(c^\prime_k) \ot a_k,$ with
$c_k,c_k^\prime \in \clc_0,$ $A_k \in \clb(\clh)$ and $a_k \in \a0.$ From this
it follows that $F^\prime$ is properly supported. 

It is easy to show the equivariance of $F^\prime.$ Indeed, $U(F^\prime \ot
1)U^*=(id \ot id \ot \phi)((id \ot
\Delta)(U \tilde{\pi(h)}F\tilde{\pi(h)}U^*))$ by using the fact that $(id \ot
\Delta)(U)=U_{12}U_{13}$ and $\Delta$ is a $\ast$-homomorphism. Now, since it
is easy to see using what we have proved in the earlier section that $(id \ot
id \ot \phi)((id \ot \Delta)(X))=(id \ot \phi)(X) \ot 1,$ for $X \in
\clm(\clb_0(\clh) \ot \cla),$ from which the equivariance of $F^\prime$
follows. 

Finally, we can verify that $\pi(c)(F-F^\prime)$ is compact for $c \in \clc_0$,
hence for all $c \in \clc,$ by very similar arguments as in \cite{val1},
adapted to our framework in a suitable way. We omit this part of the proof,
which is anyway straightforward. 
\qed

\vspace{5mm}
Let us make some more notational convention and note some simple but useful
facts. Since the von Neumann algebra generated by $\cla$ in $\clb(\clk)$ is
the direct sum of matrix algebras $\bar{\oplus}_{\al \in I} \cla_\al
=\bar{\oplus}_ \al M_{n_\al}$, where $\bar{\oplus}$ has been used to denote
the weak (or equivalently strong) operator closure of the algebraic direct
sum, it is clear that for any $X \in \clb(\clh) \ot \cla^{\prime
\prime}=\bar{\oplus}_\al \clb(\clh) \ot M_{n_\al}$, so in particular for $X \in
\clm(\clb_0(\clh) \ot \cla),$ we can write $X=\sum_{\al,i,j} X^\al_{ij} \ot
e^\al_{ij} $ as a strongly convergent sum, with $X^\al_{ij} \in \clb(\clh).$
For $\xi \in \clh,$ it is clear that  $X \xi =\sum_{\al i j} X^\al_{ij} \xi
\ot e^\al_{ij} \in \clh \ot \cla^{\prime \prime},$ where we recall that $\clh
\ot \cla^{\prime \prime}$ is the smallest Hilbert von Neumann $\cla^{\prime
\prime}$-module   generated by   the algebraic right $\cla$-module $\clh
\ot_{\rm alg} \cla.$     
 Recall from the Introduction  that $X \xi \in \clh \ot \clm(\cla),$ and  $(id
\ot \Delta)(X \xi)=(id \ot \Delta)(X)\xi \forall \xi \in \clh, X \in
\clm(\clb_0(\clh) \ot \cla).$ 

  Furthermore, for any finite subset $J \subseteq I,$ it is clear that
$X_J :=(1 \ot e_J)X =\sum_{\al \in J, ij=1,.....,n_\al} X^\al_{ij}.$  Now,
recall   that there is a bijection of the index set $I$, say $\al \mapsto
\al^\prime,$ such that $S(e_\al)=e_{\al^\prime}, S(e_{\al^\prime})=e_\al.$ In
particular, $S(e_\al+e_{\al^\prime})=e_\al+e_{\al^\prime}.$ Thus, given any
finite $J \subseteq I,$ we can enlarge $J$ suitably such that $S(e_J)=e_J$.
Since $e_J$'s are central projections, it is easy to see that $(id \ot
S)(U_J)=U^*_J=(U_J)^*$ whenever $S(e_J)=e_J.$      Now, for $\xi \in \clh_0,$
say $\xi=\pi(c) \eta$, $c \in \clc_0, \eta \in \clh,$ we have that
$(\tilde{\pi(h)}U)\xi=((\pi \ot id)((h \ot 1)\Delta_\clc(c))U)\eta=((1 \ot
e_J)(\pi \ot id)((h \ot 1)\Delta_\clc(c))U)\eta,$ where we have chosen some
finite subset $J$ of $I$ by using the fact that $(h \ot 1)\Delta_\clc(c) \in
\clc_0 \ot \a0,$ and if necesasary by enlarging $J$ suitably, assumed that
$S(e_J)=e_J.$ Thus, we can write $(\tilde{\pi(h)}U)\xi $ as a finite sum over
some set indexed by $p$ (say) of the form $\sum_p \pi(h)U_1^{(p)} \xi \ot
U^{(p)}_2,$ where $U^{(p)}_1 \in \clb(\clh), U^{(p)}_2 \in \a0,$ and we also
have $\sum_p  U_1^{(p)^*}  \ot U^{(p)^*}_2=\sum_p U_1^{(p)}  \ot
S(U^{(p)}_2).$ 

Let $\clh_0:=\pi(\clc)\clh.$ By the fact that $F^\prime$ is properly
supported, it is clear that $F^\prime \clh_0 \subseteq \clh_0.$  We now equip
$\clh_0$ with a right $\hat{\a0}$-module structure. Define  $$(\xi.a):=(id \ot
\psi_{\theta^{-1}S(a)\theta^{-2}})(U) \xi,$$  for $\xi \in \clh_0, a \in \a0.$
It is useful to note that for $c \in \clb(\clh) \ot_{\rm alg}  \a0,$ $(id \ot
\psi_{\theta^{-1}S(a)\theta^{-2}})(c)=(id \ot (\psi_{\theta
a} \circ S^{-1}))(c)$ by simple calculation using the properties of $\psi$ and
$\theta$ described in the previous section. By taking suitable limit, it is
easy to extend this for $c \in \clm(\clb_0(\clh) \ot \cla)$, in particular for
$U$. So we also have that $\xi .a =(id \ot \psi_{\theta a} \circ S^{-1})(U)
\xi.$  
\bppsn
$(\xi.a).b=\xi.(a \ast b)$ for $a,b \in \a0, \xi \in \clh_0.$ That is,
$(\xi,a) \mapsto \xi .a$ is indeed a right $\hat{\a0}$-module action.
\eppsn
 {\it Proof :-}\\
  Choosing finite subsets $J,K$ of $I$ such that $\theta^{-1}S(a) \theta^{-2}
\in supp(e_K), \theta^{-1}S(b)\theta^{-2} \in supp(e_J),$ we have that 
\bean
\lefteqn{(\xi.a).b}\\
&=& \sum_{\al \in J; i,j=1,..,n_\al}
U^\al_{ij}(\xi.a)\psi_{\theta^{-1}S(b)\theta^{-2}}(e^\al_{ij})\\ 
&=& \sum_{\al \in J;i,j=1,...,n_\al}
\sum_{\beta \in K;k,l=1,...,n_\beta} U^\al_{ij}U^\beta_{kl} \xi
\psi_{\theta^{-1}S(b)\theta^{-2}}(e^\al_{ij})\psi_{\theta^{-1}S(a)\theta^{-2}}(e^\beta_{kl})\\
&=& (id \ot \psi_{\theta^{-1}S(b)\theta^{-2}} \ot
\psi_{\theta^{-1}S(a)\theta^{-2}})(U_{12}U_{13})\xi \\
&=&  (id \ot \psi_{\theta^{-1}S(b)\theta^{-2}} \ot
\psi_{\theta^{-1}S(a)\theta^{-2}})((id \ot \Delta)(U))\xi\\
&=& (id \ot (\psi_{\theta^{-1}S(b)\theta^{-2}} \ast
\psi_{\theta^{-1}S(a)\theta^{-2}}))(U)\xi\\
&=& (id \ot \psi_{(\theta^{-1}S(b)\theta^{-2})\ast
(\theta^{-1}S(a)\theta^{-2})})(U) \xi.\\ \eean
Now, by a straightforward calculation using the properties of $\psi,$ $S$ and
$\theta$ one can verify that $(\theta^{-1}S(b)\theta^{-2})\ast
(\theta^{-1}S(a)\theta^{-2})=\theta^{-1}S(a \ast b)\theta^{-2},$ which
completes the proof. 
\qed

\vspace{5mm}
For $\xi,\eta \in \clh_0,$ say of the form
$\xi=\pi(c_1)\xi^\prime,\eta=\pi(c_2)\eta^\prime,$ it is clear that
$T_{\xi \eta}(U)$ is an element of $\a0,$ since
$(\tilde{\pi(c_1^*)}U\tilde{\pi(c_2)}=(\pi \ot id)((c_1^* \ot
1)\Delta_\clc(c_2))U,$ which belongs to $  (\pi(\clc_0) \ot_{\rm alg}
\a0)\clm(\clb_0(\clh) \ot \cla) \subseteq \clb(\clh) \ot_{\rm alg} \a0.$ We
define $$ <\xi,\eta>_{\hat{\a0}}:=\theta^{-1}T_{\xi \eta}(U) \in \hat{\a0},$$
 identifying $\a0$ as the $\ast$-algebra $\hat{\a0}$ described earlier.

We shall show that $\clh_0$ with the above right $\hat{\a0}$-action and the
$\hat{\a0}$-valued bilinear form $<.,.>_{\hat{\a0}}$ is indeed a pre-Hilbert
$\hat{\a0}$-module. However, instead of proving it directly, we shall prove it
by embedding $\clh_0$ into the free pre-Hilbert $\hat{\a0}$-module $\clf_0
:=\clh \ot_{\rm alg} \hat{\a0}$ (with the natural $\hat{\a0}$-action given by
$(\xi \ot a)b:=\xi \ot (a \ast b), \xi \in \clh, a,b \in \a0 \equiv
\hat{\a0}$), and showing that the pull back of the natural $\hat{\a0}$-valued
inner product of $\clf_0$ (which is given by $<\xi \ot a,\eta \ot
b>:=<\xi,\eta> a^\sharp \ast b,$ $\xi,\eta \in \clh, a,b\in \a0 \equiv
\hat{\a0}$) coincides with $<.,.>_{\hat{\a0}}.$ 

Define $\Sigma : \clh_0 \raro \clf_0$ by $$ \Sigma(\xi):= ((\pi(h) \ot
\theta^{-1})U )\xi,$$
 for $\xi \in \clh_0.$ Note that by writing $\xi=\pi(c)\xi^\prime$ for some
$\xi^\prime \in \clh, c \in \clc_0,$ we have that $((\pi(h) \ot 1)U) \xi=
 ((\pi \ot id)((h \ot 1)\Delta_\clc(c))U)\xi^\prime,$ and since $   (\pi \ot
id)((h \ot 1)\Delta_\clc(c))U \in \pi(\clc_0) \ot_{\rm alg} \a0,$ the range of
$\Sigma$ is clearly in $\clh \ot_{\rm alg} \a0.$ We now prove that $\Sigma$ is
in fact a module map and preserves the bilinear form $<.,.>_{\hat{\a0}}$ on
$\clh_0.$ \bppsn
\label{sigma}
For $\xi,\eta \in \clh_0, a \in \a0,$ we have that\\
(i) $\Sigma(\xi .a)=\Sigma(\xi)a.$\\
(ii) $<\Sigma(\xi),\Sigma(\eta)>=<\xi,\eta>_{\hat{\a0}}.$
 \eppsn
{\it Proof :-}\\
(i)  
Choose suitable finite set indexed by $p$ such that
$(\tilde{\pi(h)}U)\xi=\sum_p \pi(h)U_1^{(p)} \xi \ot U^{(p)}_2,$ where
$U^{(p)}_1 \in \clb(\clh), U^{(p)}_2 \in \a0,$ and  also  $\sum_p 
U_1^{(p)^*}  \ot U^{(p)^*}_2=\sum_p U_1^{(p)}  \ot S(U^{(p)}_2).$ 
Using the facts that
$\Delta(\theta)=\theta \ot \theta,$ $S^{-1}(\theta)=\theta^{-1}$ and
that $\psi(b\theta)=\psi(\theta b) \forall b \in \a0,$ and also the easily
verifiable relation $\psi_a \circ S^{-1}=\psi_{\theta^{-1}S(a)\theta^{-1}}$
for $a \in \a0,$ we have that   
\bean  
\lefteqn{ \Sigma(\xi)a}\\ 
&=&   \sum_p \pi(h)U_1^{(p)} \ot (id \ot
\psi_{\theta^{-1}S(a)\theta^{-1}})(\Delta(\theta^{-1}U_2^{(p)}))\\
&=& \sum_p \pi(h)U_1^{(p)} \ot (id \ot
\psi_{\theta^{-1}S(a)\theta^{-1}})( (\theta^{-1} \ot
\theta^{-1})\Delta(U_2^{(p)}))\\
&=& (\pi(h) \ot \theta^{-1} \ot \psi_{\theta^{-1}S(a)\theta^{-2}})((\id \ot
\Delta)(U\xi))\\
&=& (\pi(h) \ot \theta^{-1} \ot \psi_{\theta^{-1}S(a)\theta^{-2}})((id \ot
\Delta)(U) \xi)\\
&=& (\pi(h) \ot \theta^{-1} \ot
\psi_{\theta^{-1}S(a)\theta^{-2}})(U_{12}U_{13})\xi \\
&=& (\pi(h) \ot \theta^{-1})(U) ((id \ot
\psi_{\theta^{-1}S(a)\theta^{-2}})(U)\xi)\\
&=& (\pi(h) \ot \theta^{-1})(U)(\xi.a)\\
&=& \Sigma(\xi a).
\eean
(ii) Choosing suitable finite index sets as explained before, such that 
 $(\pi(h) \ot 1)U \xi=\sum_p U_1^{(p)} \ot U_2^{(p)},$ with $\sum_p
U_1^{(p)} \ot S(U_2^{(p)})=\sum_p U_1^{(p)^*} \ot U_2^{(p)^*},$ and similarly
for $(\pi(h) \ot 1)U \eta$ with the index $p$ replaced by say $q$, we can
write  
\bean 
\lefteqn{<\Sigma(\xi),\Sigma(\eta)>}\\ 
&=& \sum_{p,q}
<U_1^{(p)}\xi, \pi(h^2)U_1^{(q)}\eta>(\theta^{-1}U_2^{(p)})^\sharp \ast
(\theta^{-1}U_2^{(q)})\\
&=& \sum_{p,q} <\xi, U_1^{(p)^*}
\pi(h^2)U_1^{(q)}\eta>(\theta^{-2}S^{-1}(\theta^{-1})S^{-1}(U_2^{(p)^*}))
\ast (\theta^{-1}U_2^{(q)})\\
&=& \sum_{p,q} <\xi, U_1^{(p)^*}
\pi(h^2)U_1^{(q)}\eta>(\theta^{-1}S^{-1}(U_2^{(p)^*}))
\ast (\theta^{-1}U_2^{(q)})\\
&=& \sum_{p,q} <\xi, U_1^{(p)}\pi(h^2)U_1^{(q)}\eta>\theta^{-1}(U_2^{(p)}\ast
 U_2^{(q)})\\
&=& \sum_{p,q} <\xi, U_1^{(p)}\pi(h^2)U_1^{(q)}\eta> (\phi \ot
\theta^{-1})((U_2^{(p)}\ot 1)(S \ot id)(\Delta(U_2^{(q)})))...(1),
 \eean
using the fact that $\sum_p U_1^{(p)^*} \ot U_2^{(p)^*}=\sum_p U_1^{(p)} \ot
S(U_2^{(p)})$ and the simple   observation  that
$(\theta^{-1}x) \ast (\theta^{-1}y)=\theta^{-1}(x \ast y).$ Now, 
\bean
\lefteqn{\sum_q \pi(h^2)U_1^{(q)}\eta \ot ((S \ot id)(\Delta(U_2^{(q)})))}\\
&=& (\pi(h^2) \ot S \ot id)((id \ot \Delta)(U\eta))\\
&=& (\pi(h^2) \ot S \ot id)((id \ot \Delta)(U) \eta)\\
&=& (\pi(h^2) \ot id \ot id)((U^*)_{12}U_{13})\eta ....(2).
\eean
 Thus, from (1) and (2), $<\Sigma(\xi),\Sigma(\eta)>= 
(T_{\xi \eta} \ot \phi \ot \theta^{-1})((U(\pi(h^2) \ot 1)U^* \ot
1)U_{13})=\theta^{-1}T_{\xi \eta}(U),$ since $(id \ot \phi)(U \pi(h^2)
U^*)=1.$ This completes the proof.
\qed    

Note that from the above proposition it follows in particular that $<\xi,\eta
a>_{\hat{\a0}}=<\Sigma(\xi),\Sigma(\eta
a)>=<\Sigma(\xi),\Sigma(\eta)a>=<\Sigma(\xi),\Sigma(\eta)>\ast a
=<\xi,\eta>_{\hat{\a0}} \ast a.$ Similarly,
$<\xi,\eta>_{\hat{\a0}}^\sharp=<\eta,\xi>_{\hat{\a0}},$ and $<\xi,\xi>$ is a
nonnegative element in the $\ast$-algebra $\hat{\a0}$, since $<.,.>$ on
$\clf_0$ is a nonnegative definite form.  

 Given any $C^*$-algebra which contains $\a0$ as a dense $\ast$-subalgebra, we
can complete $\clf_0$ w.r.t. the corresponding norm to get a Hilbert
$C^*$-module in which $\clf_0$ sits as a dense submodule. Let us denote by
$\clf$ and $\clf_r$ the Hilbert $\hat{\cla}$ and $\hat{\cla}_r$-modules
respectively obtained in the above mentioned procedure, by considering $\a0$
as dense $\ast$-subalgebra of $\hat{\cla}$ and $\hat{\cla}_r$ respectively. 
The corresponding completions of $\clh_0$ will be denoted by $\cle$ and
$\cle_r$ respectively. By construction, $\Sigma$ extends to an isometry from
$\cle$ to $\clf$ and also from $\cle_r$ to $\clf_r$. We denote both
these extensions by the same notation $\Sigma$, as long as no confusion
arises. Clearly, $\cle \cong \Sigma \cle \subseteq \clf$ as closed submodule,
and similar statement will be true for $\cle_r$ and $\clf_r.$

Let us now compute the explicit form of $\Sigma^*.$ Fix $\xi,\eta \in \clh_0$
and $a \in \a0.$ Using the same notation as in the proof of the Proposition
\ref{sigma}, and using the easy observation that
$(\theta^{-1}x)^\sharp=\theta^{-1} x^\sharp$ for $x \in \a0,$ we have that 
\bean 
\lefteqn{<\Sigma(\xi),\eta \ot a>}\\
&=& \sum_p <\pi(h)U^{(p)}_1 \xi, \eta>(\theta^{-1} S^{-1}(U_2^{(p)^*}) \ast a\\
&=& \sum_p <\xi, U_1^{(p)^*} \pi(h) \eta>(\theta^{-1} S^{-1}(U_2^{(p)^*}))
\ast a \\
&=& (\theta^{-1} T_{\xi,\pi(h)\eta}(U)) \ast a ,
\eean
using the fact that $\sum_p U_1^{(p)^*} \ot S^{-1}(U_2^{(p)^*})= \sum_p
U_1^{(p)} \ot U_2^{(p)}.$ Now, $\theta^{-1} T_{\xi, \pi(h)\eta}(U) \ast
a=<\xi, \pi(h)\eta>_{\hat{\a0}} \ast a =<\xi, (\pi(h)\eta)a >_{\hat{\a0}}.$
 Thus, $$ \Sigma^*(\eta \ot a)=(\pi(h) \eta)a=(id \ot
\psi_{\theta^{-1}S(a)\theta^{-2}})(U) \pi(h) \eta.$$

Let us now prove the following important result.
\bthm
\label{compact}
Let $T \in \clb(\clh)$ be equivariant, i.e.  $U(T \ot
1)U^*=T \ot 1,$ and also assume that it satisfies the following condition
 which is slightly weaker than being properly supported :\\
For $c \in \clc_0,$
one can find $c_1,...,c_m \in \clc_0, A_1,...,A_m
\in \clb(\clh)$ (for some integer $m$) such that $T \pi(c)=\sum_k
\pi(c_k)A_k.$\\
  Then we have the following :\\ 
(i) $T(\xi a)=(T\xi)a$ $ \forall
a \in \a0,$ and thus $T$ is a module map on the $\hat{\a0}$-module $\clh_0.$
Furthermore, if $T$ is self-adjoint in the sense of Hilbert space, then
$<\xi,T\eta>_{\hat{\a0}}=<T\xi,\eta>_{\hat{\a0}}$ for $\xi,\eta \in \clh_0.$ \\
(ii) $T$ is continuous in the norms of  $\cle$ as well as $\cle_r$, thus admits
continuous extsnsions on both $\cle$ and $\cle_r$. We shall denote these
extensions by $\clt$ and $\clt_r$ respectively.\\ 
(iii) If $T \pi(h)$ is
compact in the Hilbert space sense, i.e. in $\clb_0(\clh),$ then $\clt$ and
$\clt_r$ are compact in the Hilbert module sense. 
\ethm 
{\it Proof :}\\ 
(i) is
obvious from the defintion of the right $\hat{\a0}$ action, the definition of
$<.,.>_{\hat{\a0}}$, and the equivariance of $T$. Let us prove (ii) and (iii)
only for $\clt$, as the proof for $\clt_r$ will be exactly the same. In fact,
it is enough to show that $\Sigma \clt \Sigma^*$ is continuous on $\clf$, and
is compact if $T \pi(h)$ is compact in the Hilbert space sense. Let us
introduce the following notation : for $X \in \clm(\clb_0(\clh) \ot \cla), a
\in \a0, \eta \in \clh,$ define $X \ast b :=(id \ot id \ot \psi_b \circ
S^{-1})((id \ot \Delta)(X)),$ and $X \ast ( \eta \ot a):=(X \ast a)\eta.$ Note
that clearly $X \ast a \in \clm(\clb_0(\clh) \ot \cla)$, so $(X \ast a) \eta$
makes sense. Now, we observe using the equivariance of $T$ and the explicit
formula for $\Sigma^*$ derived earlier that for $\eta \in \clh, a \in \a0,$  
\bean \lefteqn{\Sigma \clt \Sigma^*(\eta \ot a)}\\
&=& (\pi(h) \ot 1)U \beta,
\eean
where $\beta \in \clh$ is given by $\beta=(id \ot
\psi_{\theta^{-1}S(a)\theta^{-2}})(U)(T\pi(h)\eta).$  Now, by using the fact
that $(id \ot \Delta)(U)=U_{12}U_{13},$ it follows by a straightforward
computation that 
\bean
\lefteqn{ (1 \ot \theta^{-1})U \beta}\\
&=& (id \ot \theta^{-1} \ot (\psi_{\theta a} \circ S^{-1}))((id \ot
\Delta)(U) (T \pi(h) \eta).
\eean
   But $ \psi_{\theta a}( S^{-1}(b))=\psi(\theta a
S^{-1}(b))=\psi(aS^{-1}(b)\theta)=\psi(aS^{-1}(\theta^{-1}b))=(\psi_a \circ
S^{-1})(\theta^{-1}b),$ and hence $(id \ot \theta^{-1} \ot (\psi_{\theta a}
\circ S^{-1} ))((id \ot \Delta)(U))=(id \ot id \ot (\psi_a \circ S^{-1}))((id
\ot \Delta)((1 \ot \theta^{-1})U))=((1 \ot \theta^{-1})U) \ast a.$ From this,
it is clear that  
$$ (\Sigma \clt \Sigma^*)(\eta \ot a)=((\pi(h) \ot \theta^{-1})U(T \pi(h) \ot
1))\ast (\eta \ot a).$$
 Now, note that $T \pi(h)=\sum_{k=1}^m
\pi(c_k)A_k,$ for some $c_1,...,c_m \in \clc_0, A_1,...,A_m \in \clb(\clh),$
and so we have  $(\pi(h) \ot \theta^{-1})U(T \pi(h) \ot
1) =\sum_k (1 \ot \theta^{-1})(\pi \ot
id)((h \ot 1)\Delta_\clc(c_k))U(A_k  \ot 1).$  But 
$(h \ot 1)\Delta_\clc(c_k)$ is in $\clc_0 \ot_{\rm alg} \a0$ for each
$k=1,..m,$  and thus   $(\pi(h) \ot \theta^{-1})U(T \pi(h) \ot
1) \in \clb(\clh) \ot_{\rm alg} \a0$ clearly. Choosing some large enough
 finite subset  $J$ of $I$ such that $(\pi(h) \ot \theta^{-1})U(T \pi(h) \ot
1) =\sum_{\al \in J, ij=1,..,n_\al} B^\al_{ij} \ot e^\al_{ij}$, (with
$B^\al_{ij} \in \clb(\clh)$), it is clear that $\Sigma \clt \Sigma^*=\sum_{\al
\in J, ij=1,...,n_\al} B^\al_{ij} \ot L_{e^\al_{ij}}$, where for $x \in \a0,$
$L_x : \a0 \raro \a0$ with $L_x(a)=x \ast a.$ As $L_x$ is a norm-continuous
map on $\hat{\cla},$ the above finite sum shows that $\Sigma \clt \Sigma^*$
indeed admits a continuous extension on the Hilbert $\hat{\cla}$-module
$\clf$. This proves (ii). 

Furthermore, since $\clk(\clh \ot \hat{\cla}) \cong \clb_0(\clh) \ot
\hat{\cla}$, where $\clk(E)$ means the set of compact (in the Hilbert module
sense) opeartors on the Hilbert module $E$, it is easy to see that $\Sigma
\clt \Sigma^*$ is compact on $\clf$ if $B^\al_{ij}$'s are compact on the
Hilbert space $\clh.$ Now, $B^\al_{ij}=(id \ot \phi^\al_{ij})
((\pi(h) \ot \theta^{-1})U(T \pi(h) \ot 1))=\pi(h).(id \ot \phi^\al_{ij})((1
\ot \theta^{-1})U)T \pi(h),$ where $\phi^\al_{ij} $ is the functional on $\a0$
which is $0$ on all $e^\beta_{kl}$ except $\beta=\al, (kl)=(ij),$  with
$\phi^\al_{ij}(e^\al_{ij})=1.$ It follows that $B^\al_{ij}$ 's are all compact
if $T \pi(h)$ is so, which completes the proof. 
\qed

Now, let us come to the construction of the Baum-Connes maps $\mu_1 :
KK_1^\cla(\clc,\IC) \raro KK_1(\IC,\hat{\cla})$ and $\mu_1^r :
KK_1^\cla(\clc,\IC) \raro KK_1(\IC,\hat{\cla}_r).$ Let us do it only for
$\mu_1$, as the case of $\mu^r_1$ is similar, and in fact $\mu^r_1$ will be
the compositon of $\mu_1$ and the canonical map from $KK_1(\IC , \hat{\cla})$
to $KK_1(\IC,\hat{\cla}_r)$ induced by the canonical surjective
$C^*$-homomorphism from $\hat{\cla}$ to $\hat{\cla}_r.$ 
Note that an element of $KK_1(\IC,\hat{\cla}) \cong K_1(\hat{\cla})$  is given
by the suitable homotopy class $[E,L]$ of a pair of the form $(E, L)$, where
$E$ is a Hilbert $\hat{\cla}$-module and $L \in \cll(E)$ (the set of
adjointable $\hat{\cla}$-linear maps on $E$) such that $L^*=L,$ $L^2-1$ is
compact in the sense of Hilbert module. For more details, see for example
\cite{kk}.

\bthm
\label{defn}
Given a cycle $(U, \pi, F) \in KK_1^\cla(\clc,\IC),$ let $F^\prime \equiv
F^\prime_h$ be the equivariant and properly supported operator as constructed
in \ref{proper}, with a given choice of $h$ as in that theorem. Then the
continuous extension of $F^\prime_h$ on the Hilbert module $\cle$ (as
described by  the Theorem \ref{compact}), to be denoted  by say
$\clf^\prime_h$, satisfies the conditions that
$(\clf^\prime_h)^*=\clf^\prime_h$ (as module map), and  $(\clf^\prime_h)^2-I$
is compact on $\cle$. Define  $$ \mu_1((U,\pi,F)):=[\cle, \clf^\prime_h] \in
KK_1(\IC,\hat{\cla}) \cong K_1(\hat{\cla}).$$  
In fact, $[\cle,\clf^\prime_h]$  is independent (upto operatorial homotopy) of
the choice of $h$.  
\ethm 
{\it Proof :-}\\
  Since $F^\prime_h$ is equivariant and properly supported, it is clear that
$T_h:=(F^\prime_h)^2-1$ is equivariant and for any $c \in \clc_0,$ there are
finitely many $c_1,...,c_m \in \clc_0, A_1,...,A_m \in \clb(\clh)$ such that
$T_h \pi(c)=\sum_k \pi(c_k)A_k.$ Furthermore, by the Theorem \ref{proper}, we
have that $\pi(c)T_h$, and hence $T_h \pi(c)$ is compact operator on $\clh$
 for every $c \in \clc.$ So, in particular, $T_h \pi(h)$ is compact. By
Theorem \ref{compact}, it follows that the continuous extension of $T_h$ on
$\cle$ is compact in the sense of Hilbert modules. Furthermore, the fact
that $(\clf^\prime_h)^*=\clf^\prime_h$ is clear from (i) of the Theorem
\ref{compact}. So, $[\cle,\clf^\prime_h] \in KK_1(\IC,\hat{\cla}).$
Furthermore, as we can see from the proof of   the Theorem \ref{proper},
$(F^\prime_h-F)\pi(c) \in \clb_0(\clh)$ $\forall c \in \clc_0,$ and so for
$h,h^\prime$ satisfying {\bf A3}, we have
$(F^\prime_h-F^\prime_{h^\prime})\pi(c) \in \clb_0(\clh),$ and hence by
Theorem \ref{compact}, $\clf^\prime_h-\clf^\prime_{h^\prime}$ is compact in
the Hilbert module sense. Thus, for each $t \in [0,1],$  setting
$\clf(t):=t\clf^\prime_{h^\prime}+(1-t)\clf^\prime_h,$ we have that   
$\clf(t)^2-I $ is compact on $\cle$, and this gives a homotopy in
$KK_1(\IC,\hat{\cla})$ between $[\cle,\clf^\prime_h]$ and
$[\cle,\clf^\prime_{h^\prime}].$  
\qed        

\vspace{5mm}

{\it Acknowledgement :}\\
 D. Goswami  would like to express his gratitude to I.C.T.P. for a visting
research fellowship during January-August 2002, and to  A.O. Kuku and  the
other organisers of the ``School and Conference on  Algebraic K Theory and Its
Applications" at I.C.T.P. (Trieste) in July 2002. He would also like to thank
  T. Schick  for sending some relevant preprint,  and   A.
Valette and  I. Chatterjee for giving useful information regarding
 some  manuscript (yet to be published)  by A. Valette.

\end{document}